\documentclass[12pt]{amsart}

\usepackage{enumerate}
\usepackage{fullpage}
\usepackage{hyperref}
\usepackage[frozencache]{minted}
\usepackage{textgreek}

\usepackage{mdframed}
\surroundwithmdframed{minted}

\mdfsetup{skipabove=1em,skipbelow=1em}

\newcommand{\RR}{\mathbb R}
\newcommand{\ZZ}{\mathbb Z}

\DeclareMathOperator{\supp}{supp}

\theoremstyle{definition}

\newtheorem{example}{Example}[section]

\title{Probability package for Macaulay2}
\author{Douglas A. Torrance}
\address{Piedmont University, Department of Mathematical Sciences, PO Box 10, Demorest, GA, 30535}
\email{dtorrance@piedmont.edu}

\begin{document}
\maketitle

\begin{abstract}
  We introduce the \verb|Probability| package for Macaulay2, which provides an interface for users to compute probabilities and generate random variates from a wide variety of univariate probability distributions.
\end{abstract}

\section{Introduction}
Users of Macaulay2 \cite{M2} often have reason to use random elements.  Entire packages, e.g., \verb|SpaceCurves| \cite{zhang}, exist to create random algebro-geometric objects.  Also, a growing number of packages, e.g., \verb|GraphicalModels| \cite{gps}, are devoted to the exciting field of algebraic statistics.

Macaulay2's built-in \verb|random| method wraps around two functions in shared libraries that it is linked against: \verb|mpz_urandomm| from GMP \cite{gmp} for generating random integers and \verb|mpfr_urandomb| from MPFR \cite{mpfr} for generating random real numbers.  Both of these use the Mersenne Twister pseudorandom number generator \cite{matsumoto-nishimura} for sampling random variates from uniform discrete and continuous probability distributions, respectively.

There has previously been no support for other probability distributions, and users interested in these distributions would have needed to look elsewhere, for example to the statistical software R \cite{R}

However, beginning with version 1.20, Macaulay2 has been distributed with the \verb|Probability| package, written by the author, that implements many common probability distributions and provides an interface for users to create additional distributions for use in their work.

This paper is outlined as follows.  In Section \ref{basics}, we review the basics of probability theory.  In Section \ref{probability-distribution class}, we introduce the \texttt{ProbabilityDistribution} class from the package and its methods.  Finally, in Section \ref{chi-squared example}, we use the package to perform an example of Pearson's chi-squared test for independence.

The author would like to thank the anonymous referees for helpful comments on both this paper and the package.

\section{Basics of probability theory}
\label{basics}

The basics of probability theory are well-covered in a large variety of undergraduate- and graduate-level textbooks, e.g., \cite{billingsley,mendenhall-scheaffer,ross}, but we review the essentials here.

Suppose $\Omega$ is a \textit{sample space} containing the possible outcomes of some random experiment, $\mathcal F$ is a $\sigma$-algebra on $\Omega$ (known as the \textit{event space}), and $P:\mathcal F\to[0,1]$ is a $\sigma$-additive \textit{probability measure} (or \textit{probability distribution}) satisfying $P(\emptyset)=0$ and $P(\Omega)=1$.  The triple $(\Omega,\mathcal F,P)$ is known as a \textit{probability space}.  A \textit{random variable} is measurable function $X:\Omega\to\RR$.  The probability distribution of this random variable is the pushforward measure $X_*P$ on $\RR$, i.e., for any Borel set $A\subset\RR$, $X_*P(A) = P(\{\omega:X(\omega)\in A)\})$, usually denoted $P(X\in A)$.  The smallest closed $A$ for which $P(X\in A)=1$ is the \textit{support} of $X$, denoted $\supp(X)$.  The Radon-Nikodym derivative of $X_*P$ with respect to the measure $\mu$ on $\RR$, i.e., the function $f_X$ satisfying $P(X\in A)=\int_Af_X\,d\mu$, is known as the \textit{probability density function} of $X$.  The function $F_X$ defined by $F_X(x) = P(X\leq x)$ is the \textit{cumulative distribution function} of $X$.

There are two main cases.  When $X(\Omega)$ is countable and $\mu$ is the counting measure, i.e., $\int_Af_X\,d\mu = \sum_{x\in{A}}f_X(x)$,  then $X$ is said to have a \textit{discrete} distribution.  In this case, $f_X$ is often referred to instead as the \textit{probability mass function}.  When $X(\Omega)$ is uncountable and $\mu$ is the usual Lebesgue measure, then $X$ has a \textit{continuous} distribution.  See Tables \ref{discrete distributions} and \ref{continuous distributions} for common discrete and continuous probability distributions, respectively.

\begin{table}
  \centering
  \begin{tabular}{|c|c|c|c|}
    \hline
    distribution & parameters & $\supp(X)$ & $f_X(x)$, $x\in\supp(X)$ \\
    \hline\hline
    binomial & $n, p$ & $\{0,\ldots,n\}$ & $\dbinom{n}{x}p^x(1-p)^{n-x}$ \\[1em]
    Poisson & $\lambda $ & $\{0,1,\ldots\}$ & $\dfrac{\lambda^x}{x!}e^{-\lambda}$ \\[1em]
    geometric & $p$ & $\{0,1,\ldots\}$ & $p(1-p)^x$ \\[1em]
    negative binomial & $r, p$ & $\{0,1,\ldots\}$ & $\dfrac{\Gamma(x + r)}{\Gamma(r) x!} p^r (1 - p)^x$ \\[1em]
    hypergeometric & $m,n,k$ & $\{0,\ldots,m\}$ & $\dfrac{\binom{m}{x}\binom{n}{k-x}}{\binom{m+n}{k}}$ \\[1em]
    \hline
    \end{tabular}
    \caption{Common discrete probability distributions}
  \label{discrete distributions}
\end{table}

\begin{table}
  \centering
  \begin{tabular}{|c|c|c|c|}
    \hline
    distribution & parameters & $\supp(X)$ & $f_X(x)$, $x\in\supp(X)$ \\
    \hline\hline
    uniform & $a, b$ & $[a,b]$ & $\dfrac{1}{b-a}$ \\[1em]
    exponential & $\lambda$ & $[0,\infty)$ & $\lambda e^{-\lambda x}$ \\[1em]
    normal & $\mu,\sigma$ & $\RR$ & $\dfrac{1}{\sigma\sqrt{2\pi}}e^{-(x-\mu)^2/(2\sigma^2)}$ \\[1em]
    gamma & $\alpha,\lambda$ & $[0,\infty)$ & $\dfrac{\lambda^\alpha}{\Gamma(\alpha)}x^{\alpha - 1}e^{-\lambda x}$ \\[1em]
    chi-squared & $\nu$ & $[0,\infty)$ & $\dfrac{1}{2^{\nu/2}\Gamma(\nu/2)} x^{\nu/2 - 1}e^{-x/2}$ \\[1em]
    Student's $t$ & $\nu$ & $\RR$ &                             $\dfrac{\Gamma((\nu+1)/2)}{\sqrt{\nu\pi}\Gamma(\nu/2)}(1 + x^2/\nu)^{-(\nu + 1)/2}$ \\[1em]
    $F$ & $\nu_1,\nu_2$ & $[0,\infty)$ & $\sqrt{\dfrac{(\nu_1x)^{\nu_1}\nu_2^{\nu_2}}{(\nu_1x+\nu_2)^{\nu_1+\nu_2}}}/(xB(\frac{\nu_1}{2}, \frac{\nu_2}{2}))$ \\[1em]
    beta & $\alpha$, $\beta$ & $[0, 1]$ & $\dfrac{x^{\alpha - 1}(1 - x)^{\beta - 1}}{B(\alpha, \beta)}$ \\[1em]
    \hline
    \end{tabular}
    \caption{Common continuous probability distributions}
    \label{continuous distributions}
\end{table}

The \textit{quantile function} is the function $Q_X:[0,1]\to\RR\cup\{-\infty,\infty\}$ defined by $Q_X(p) = \inf\{x:F_X(x)\geq p\}$.  It is useful for generating random variates of $X$ using the \textit{inversion method}.  In particular, if $u$ is a random variate of a random variable $U$ with the uniform distribution on the interval $[0, 1]$, then $Q_X(u)$ is a random variate of $X$ \cite{devroye}.

\section{The \texttt{ProbabilityDistribution} class}
\label{probability-distribution class}

The \texttt{Probability} package defines an abstract \texttt{ProbabilityDistribution} class with two subclasses, \texttt{DiscreteProbabilityDistribution} and \texttt{ContinuousProbabilityDistribution}.  Discrete and continuous random variables are represented by instances of these two classes, respectively.

There are four installed methods.  If \texttt{X} is a \texttt{ProbabilityDistribution} object corresponding to a random variable $X$, then

\begin{enumerate}[(a)]
\item \verb|density_X| is the probability density function $f_X$,
\item \verb|probability_X| is the cumulative distribution function $F_X$,
\item \verb|quantile_X| is the quantile function $Q_X$, and
\item \verb|random X| is a random variate of $X$.
\end{enumerate}

To construct a \texttt{ProbabilityDistribution} object, use of one of the constructor methods, \texttt{discreteProbabilityDistribution} and \texttt{continuousProbabilityDistribution}.

At minimum, these methods take a function (the probability density function) as input.  If the support differs from $\{0, 1, \ldots\}$ (for discrete probability distributions) or $[0, \infty)$ (for continuous probability distributions), then it may be specified using the \verb|Support| option.  This option takes a sequence of the form $(a, b)$.  In the discrete case, this is interpreted as the set $\{a,\ldots,b\}\subset\ZZ$, and in the continuous case, the interval $(a,b)\subset\RR$.  Setting $a=-\infty$ and/or $b=\infty$ is also allowed.

A cumulative distribution function based on the probability density function is installed by default (using \verb|sum| for discrete probability distributions and \verb|integrate| for continuous probability distributions), but one may also be specified using the \verb|DistributionFunction| option.

Similarly, a na\"ive quantile function (adding probabilities inside a \verb|while| loop until the target probability $p$ is obtained in the discrete case and using the bisection method for numerically approximating a root of $F_X(x) - p$ in the continuous case) is installed by default, but may be specified using the \verb|QuantileFunction| option.

By default, \verb|random| uses the inversion method to generate random variates, but this may be overriden using the \verb|RandomGeneration| option.

\begin{example}[Six-sided die]
  Consider the following simple example, where $X$ is a discrete random variable corresponding to the result of rolling a fair six-sided die,  demonstrating each of the above methods.

\begin{minted}{macaulay2}
i1 : needsPackage "Probability";

i2 : X = discreteProbabilityDistribution(x -> 1/6, Support => (1, 6));

i3 : density_X 3

     1
o3 = -
     6

o3 : QQ

i4 : probability_X 3

     1
o4 = -
     2

o4 : QQ

i5 : quantile_X oo

o5 = 3

i6 : random X

o6 = 2
\end{minted}
\end{example}

\begin{example}[Triangular distribution]
  We next give an example of a continuous probability distribution, the Bates distribution of the mean of two random variables with the uniform distribution on the unit interval, which has probability density function
  \begin{equation*}
    f_X(x) = \begin{cases}
      4x &\text{if }0\leq x\leq\frac{1}{2} \\
      4(1 - x) &\text{if }\frac{1}{2} <x\leq 1\\
        0 &\text{otherwise,}
      \end{cases}
    \end{equation*}
and so its probabilities are areas inside an isoceles triangle.

  \begin{minted}{macaulay2}
i1 : needsPackage "Probability";

i2 : X = continuousProbabilityDistribution(
         x -> if x < 1/2 then 4*x else 4*(1 - x),
         Support => (0, 1));

i3 : density_X(2/3)

     4
o3 = -
     3

o3 : QQ

i4 : probability_X(2/3)

o4 = .777777777777778

o4 : RR (of precision 53)

i5 : quantile_X oo

o5 = .666666666666668

o5 : RR (of precision 53)

i6 : random X

o6 = .318199004810108

o6 : RR (of precision 53)
\end{minted}

\end{example}

As is the case when working over $\RR$ in any computer system, there can be unexpected results due to floating point rounding errors.  This is apparent in the outputs \verb|o4| and \verb|o5| in the example above, where outputs of $\frac{7}{9}$ and $\frac{2}{3}$ would be expected.  Users should take care to be aware of this.

Alternatively, one of the built-in constructor methods for each of the common probability distributions listed in Tables \ref{discrete distributions} and \ref{continuous distributions} may be used.  For example, \verb|binomialDistribution(n, p)| will return a \verb|DiscreteProbabilityDistribution| object corresponding to a binomially distributed random variable with parameters $n$ and $p$.

A number of new special functions (e.g, \verb|Beta| and \verb|inverseErf|) were added to Macaulay2 beginning in version 1.20 by the author and Paul Zinn-Justin, wrapping around corresponding functions from the Boost Math Toolkit \cite{boost} and MPFR \cite{mpfr}.  Many of these were used to define the density, distribution, and/or quantile functions of these common distributions.

\section{Chi-squared test for independence}
In one of the seminal works in algebraic statistics, Diaconis and Sturmfels \cite{ds}  considered a data set from \cite{ah}.  A $12\times 12$ contingency table, represented by the matrix \verb|birthDeath| below, contains birth and death month information for 82 descendants of Queen Victoria.  In particular, the rows correspond to birth months (January to December) and the columns correspond to death months.  So for example, the \verb|2| in \verb|birthDeath_(3, 2)| means that two of the 82 people represented in the data set were in born in April and died in March.
\label{chi-squared example}

\begin{minted}{macaulay2}
i1 : birthDeath = matrix {
         {1, 0, 0, 0, 1, 2, 0, 0, 1, 0, 1, 0},
         {1, 0, 0, 1, 0, 0, 0, 0, 0, 1, 0, 2},
         {1, 0, 0, 0, 2, 1, 0, 0 ,0 ,0 ,0, 1},
         {3, 0, 2, 0, 0, 0, 1, 0, 1, 3, 1, 1},
         {2, 1, 1, 1, 1, 1, 1, 1, 1, 1, 1, 0},
         {2, 0, 0, 0, 1, 0, 0, 0, 0, 0, 0, 0},
         {2, 0, 2, 1, 0, 0, 0, 0, 1, 1, 1, 2},
         {0, 0, 0, 3, 0, 0, 1, 0, 0, 1, 0, 2},
         {0, 0, 0, 1, 1, 0, 0, 0, 0, 0, 1, 0},
         {1, 1, 0, 2, 0, 0, 1, 0, 0, 1, 1, 0},
         {0, 1, 1, 1, 2, 0, 0, 2, 0, 1, 1, 0},
         {0, 1, 1, 0, 0, 0, 1, 0, 0, 0, 0, 0}};

              12        12
o1 : Matrix ZZ   <--- ZZ
\end{minted}

The question at hand is whether these two variables (birth month and death month) are independent.  A standard test is due to Pearson \cite{pearson}.  In particular, given an $m\times n$ contingency table $O$, let $O_{i\bullet}=\sum_{j=1}^n O_{ij}$, $O_{\bullet j} = \sum_{i=1}^m O_{ij}$, and $N = \sum_{i=1}^m\sum_{j=1}^n O_{ij}$.  Then, defining $E_{ij} = \frac{O_{i\bullet}O_{\bullet j}}{N}$ for all $i,j$, the test statistic
\begin{equation*}
  \chi^2 = \sum_{i=1}^m\sum_{j=1}^n\frac{(O_{ij} - E_{ij})^2}{E_{ij}}
\end{equation*}
converges in distribution to a chi-squared distribution with $(m-1)(n-1)$ degrees of freedom.

First, we compute the test statistic, noting that that $O_{i\bullet}O_{\bullet j}$ and $N$ may be computed nicely by multiplying $O$ by various matrices of ones.

\begin{minted}{macaulay2}
i2 : ones = (m, n) -> matrix toList(m : toList(n : 1));

i3 : N = (ones(1, 12) * birthDeath * ones(12, 1))_(0, 0);

i4 : expected = (1/N) * birthDeath * ones(12, 12) * birthDeath;

              12        12
o4 : Matrix QQ   <--- QQ

i5 : chi2 = numeric sum flatten table(12, 12, (i, j) ->
         (birthDeath_(i, j) - expected_(i, j))^2 / expected_(i, j))

o5 = 115.559632730585

o5 : RR (of precision 53)
\end{minted}

Finally, we use the \verb|Probability| package to compute the $p$-value of the hypothesis test.

\begin{minted}{macaulay2}
i6 : needsPackage "Probability";

i7 : X = chiSquaredDistribution((12 - 1) * (12 - 1))

o7 = chi2(121)

o7 : ContinuousProbabilityDistribution

i8 : probability_X(chi2, LowerTail => false)

o8 = .622505910459144

o8 : RR (of precision 53)
\end{minted}

In particular, if $X$ is a continuous random variable with the chi-squared distribution with 121 degrees of freedom and $\chi^2$ is the Pearson test statistic computed above, then $P(X > \chi^2) = 0.6225$.  So under the null hypothesis that birth and death months are independent, it is quite likely that the sample data may have been obtained, and thus there is not sufficient evidence to reject it in favor of the alternate hypothesis that they are dependent.

\bibliography{probability-m2}{}
\bibliographystyle{abbrv}

\end{document}